\documentstyle[12pt, amsfonts]{amsart}
\def\H{{\mathbb H}}
\def\Cnh{{\C^{n,1}}}
\def\proj{{\mathbb{P}}}
\def\chs{{\H_\C^n}}
\def\C{{\mathbb{C}}}
\def\Cn{{\C^n}}
\def\Rn{{\R^n}}
\def\N{{\mathbb{N}}}
\def\s{{\mathrm{S}}}
\def\sn{{\s^{n-1}}}
\def\B{{\mathrm{B}}}
\def\Vol{{\mathrm{Vol}}}
\def\HVol{{\mathrm{HVol}}}

\def\sw{{\mathcal S}}
\def\RT{{\mathcal R}}

\def\note#1{\ifvmode\leavevmode\fi\vadjust{\vbox to0pt{\vss
 \hbox to 0pt{\hskip\hsize\hskip1em
\vbox{\hsize2cm\small\raggedright\pretolerance10000
 \noindent #1\hfill}\hss}\vbox to8pt{\vfil}\vss}}}
 
\renewcommand{\:}{\, : \,}

\begin{document}
\def\R{{\mathbb R}}
\def\Z{{\mathbb Z}}
\def\C{{\mathbb C}}
\newcommand{\trace}{\rm trace}
\newcommand{\Ex}{{\mathbb{E}}}
\newcommand{\Prob}{{\mathbb{P}}}
\newcommand{\E}{{\cal E}}
\newcommand{\F}{{\cal F}}
\newtheorem{df}{Definition}
\newtheorem{theorem}{Theorem}
\newtheorem{lemma}{Lemma}
\newtheorem{pr}{Proposition}
\newtheorem{co}{Corollary}
\def\sign{\mbox{ sign }}
\def\a{\alpha}
\def\N{{\mathbb N}}
\def\A{{\cal A}}
\def\L{{\cal L}}
\def\X{{\cal X}}
\def\F{{\cal F}}
\def\c{\bar{c}}
\def\diam{\mbox{\rm dim}}
\def\vol{\mbox{\rm Vol}}  
\def\b{\beta}
\def\t{\theta}
\def\l{\lambda}
\def\e{\varepsilon}
\def\colon{{:}\;}
\def\pf{\noindent {\bf Proof :  \  }}
\def\endpf{ \begin{flushright}
$ \Box $ \\
\end{flushright}}

\title{The Busemann-Petty Problem in Complex Hyperbolic Space}
\author{Susanna Dann}
\address{Mathematics Department\\
University of Missouri\\
Columbia, Missouri}
\email{danns@@missouri.edu}
\date{July 30, 2012}
\begin{abstract}
The Busemann-Petty problem asks whether origin-symmetric convex bodies in $\R^n$ with smaller central hyperplane sections necessarily have smaller volume. The answer is affirmative if $n\leq 4$ and negative if $n\geq 5$. We study this problem in the complex hyperbolic $n$-space $\chs$ and prove that the answer is affirmative for $n\leq 2$ and negative for $n\geq 3$.
\end{abstract}
\maketitle
\section*{Introduction}
\noindent
The Busemann-Petty problem asks the following question. Given two origin symmetric convex bodies $K$ and $L$ in $\Rn$ such that 
$$ \Vol_{n-1}(K\cap H) \leq \Vol_{n-1}(L\cap H) $$
for every hyperplane $H$ in $\Rn$ containing the origin, does it follow that 
$$ \Vol_n(K) \leq \Vol_n(L)? $$
The answer is affirmative if $n\leq 4$ and negative if $n\geq 5$. The problem, posed in 1956 in \cite{BusemannPetty1956}, was solved in the late 90's as a result of a sequence of papers \cite{LR,Ba,Gi,Bo,Lu,Pa,Ga1,Ga2,Zh1,K1,K2,Zh2, GardnerKoldobskySchlumprecht1999}, see \cite{Koldobsky2005}, p. 3-5, for the history of the solution. 

Since then the Busemann-Petty problem was studied on other spaces as were its numerous generalizations. V. Yaskin studied the Busemann-Petty problem in real hyperbolic and spherical spaces, \cite{yaskin2006}. He showed that for the spherical space the answer is the same as for $\Rn$, but not so for the real hyperbolic space, where the answer is affirmative for $n\leq 2$ and negative for $n\geq 3$. A. Koldobsky, H. K{\"o}nig and M. Zymonopoulou demonstrated in \cite{KoldobskyKonigZymonopoulou2008} that the answer to the complex version of the Busemann-Petty problem is affirmative for the complex dimension $n\leq 3$ and negative for $n\geq 4$.  Other results on the complex Busemann-Petty problem include \cite{Zymonopoulou2008,Zymonopoulou2009,rubin2010,koldobsky2011}.

In this article we consider the Busemann-Petty problem in the complex hyperbolic $n$-space. For $\xi \in \C^n$ with $|\xi|=1$, denote by 
$$ H_{\xi} := \{z\in \Cn \: (z,\xi)=\sum\limits_{k=1}^{n}z_k \overline{\xi_k} = 0\} $$
the complex hyperplane through the origin perpendicular to $\xi$. We identify $\C^n$ with $\R^{2n}$ via the mapping 
\begin{equation}\label{eq_identification}
 (\xi_{11}+i\xi_{12}, \dots, \xi_{n1}+i\xi_{n2}) \mapsto (\xi_{11}, \xi_{12}, \dots, \xi_{n1}, \xi_{n2}) \, .
\end{equation}
Under this mapping the hyperplane $H_{\xi}$ turns into a $(2n-2)$-dimensional subspace of $\R^{2n}$ orthogonal to the vectors
$$ \xi=(\xi_{11}, \xi_{12}, \dots, \xi_{n1}, \xi_{n2}) \text{ and } \xi_{\perp}=(-\xi_{12}, \xi_{11}, \dots, -\xi_{n2}, \xi_{n1}) \, . $$
A convex body $K$ in $\R^{2n}$ is called \textit{$R_{\theta}$-invariant}, if for every $\theta \in [0,2\pi]$ and every $\xi=(\xi_{11}, \xi_{12}, \dots, \xi_{n1}, \xi_{n2})\in \R^{2n}$
$$ \|\xi \|_K = \| R_{\theta}(\xi_{11}, \xi_{12}), \dots, R_{\theta}(\xi_{n1},  \xi_{n2}) \|_K \, ,$$
where $R_{\theta}$ stands for the counterclockwise rotation by an angle $\theta$ around the origin in $\R^2$.

Recall that origin symmetric convex bodies in $\Cn$ are unit balls of norms on $\Cn$ and therefore, under the mapping (\ref{eq_identification}), they are $R_{\theta}$-invariant convex bodies in $\R^{2n}$. We will work with the ball model of the $\chs$ and consequently we will only consider bodies contained in the unit ball.
We denote the volume element on $\chs$ by $d\mu_n$ and the volume of a body $K$ in $\R^{2n}$ with respect to this volume element by $\HVol_{2n}(K)$ to distinguish from the Euclidean volume of $K$.

Now the Busemann-Petty problem in $\chs$ can be posed as follows. Given two $R_{\theta}$-invariant convex bodies $K$ and $L$ in $\R^{2n}$ contained in the unit ball such that 
$$ \HVol_{2n-2}(K\cap H_{\xi}) \leq \HVol_{2n-2}(L\cap H_{\xi}) $$
for $\xi$ an element of the unit sphere $\s^{2n-1}$ of $\R^{2n}$, does it follow that 
$$ \HVol_{2n}(K) \leq \HVol_{2n}(L)? $$
Analytic solutions of the Busemann-Petty problem in different settings are based on establishing a connection between a certain distribution and the problem. E. Lutwak introduced a class of \textit{intersection bodies} in \cite{Lu} and established a connection between this class and the Busemann-Petty problem on $\Rn$. Recall that an origin symmetric star body $K$ in $\Rn$ is an intersection body if and only if $\|\cdot \|_K^{-1}$ is a positive definite distribution on $\Rn$, see Theorem 4.1 in \cite{Koldobsky2005}. Later, A. Zvavitch solved the Busemann-Petty problem on $\Rn$ for arbitrary measures, \cite{Zvavitch2005}. He linked the problem to the distribution $\|x\|_K^{-1} \frac{f_n(x \|x\|_K^{-1})}{f_{n-1}(x \|x\|_K^{-1})}$, where $f_n$, a locally integrable function on $\Rn$, is the density function for a measure $\mu_n$ on $\Rn$ and $f_{n-1}$, a function on $\Rn$ integrable on central hyperplanes, is the density function for a measure $\mu_{n-1}$ on $\R^{n-1}$. In \cite{yaskin2006} V. Yaskin established a connection between the Busemann-Petty problem in hyperbolic and spherical spaces and the distributions $\frac{\|x\|_K^{-1}}{1\pm (|x| \|x\|_K^{-1})^2}$ as a special case of Zva\-vitch's theorem. Recently, A. Koldobsky at al. found a connection between the $2$-intersection bodies and the Busemann-Petty problem on $\Cn$, \cite{KoldobskyKonigZymonopoulou2008}. The classes of $k$-intersection bodies were introduced in \cite{Koldobsky1999, Koldobsky2000}. Recall that an origin symmetric star body $K$ in $\Rn$ is a $k$-intersection body, $0<k<n$, if and only if $\|\cdot \|_K^{-k}$ is a positive definite distribution on $\Rn$, see \cite{Koldobsky2000}. This is also the case for the Busemann-Petty problem on $\chs$, namely we prove in Theorem \ref{thm_CbeAtoBRPandD} that the answer to the problem is affirmative if and only if for every $R_{\theta}$-invariant convex body in $\R^{2n}$ contained in the unit ball the distribution $\frac{\|x\|_K^{{-2}}}{1-\left(|x| \|x\|_K^{-1}\right)^2}$ is positive definite. Then we prove that the latter is true for $n=2$ and false for $n\geq 3$, providing a solution to the problem. In our proof we use methods from \cite{yaskin2006}, \cite{KoldobskyKonigZymonopoulou2008} as well as recently obtained results for complex star bodies from \cite{KoldobskyPaourisZymonopoulou2011}.  

A few generalizations of the complex Busemann-Petty problem have been considered so far, see \cite{Zymonopoulou2008, Zymonopoulou2009}. For another result in the complex hyperbolic space related to convex geometry see \cite{AbardiaGallego2011}. Some other extensions of results from convex geometry to non-Euclidean settings include \cite{GaoHugSchneider2001, Gardner2002, BezdekSchneider2010, AbardiaBernig2011}. For other generalizations of the Busemann-Petty problem see \cite{BourgainZhang1999, Koldobsky1999, Koldobsky2000, Koldobsky2002, Koldobsky2003, Zvavitch2005, KoldobskyYaskinYaskina2006}.

\section{Preliminaries}\label{SectionPreliminaries}
\subsection{Complex Hyperbolic Space}
The material of sections \ref{subsubsectionTheBallModel} and \ref{subsubsectionTheBergmanMetric} is taken from the book by Goldman \cite{goldman99}. We refer the interested reader to this book for more information. 
\subsubsection{The Ball Model}\label{subsubsectionTheBallModel}
Let $V$ be a complex vector space. The \textit{projective space associated to $V$} is the space $\proj(V)$ of all lines in $V$, i.e. one dimensional complex linear subspaces through the origin. 

Let $\C^{n,1}$ be the $(n+1)$-dimensional complex vector space consisting of $(n+1)$-tuples 
$$ Z = \left[ \begin{array}{l} 
									Z' \\
             			Z_{n+1}         
             	\end{array} \right] \in \C^{n+1}	$$ 
and equipped with the indefinite \footnote{neither positive- nor negative-semidefinite} Hermitian form 
\begin{align*}
     \left\langle Z, W \right\rangle 	&:= (Z', W')- Z_{n+1} \overline{W}_{n+1} \\
     																	&= Z_1 \overline{W}_1 + \dots + Z_n \overline{W}_n - Z_{n+1} \overline{W}_{n+1} \, ,
\end{align*}
where $Z'$ is a vector in $\C^n$ and $Z_{n+1}\in\C$. Consider the subset of \textit{negative vectors of $\C^{n,1}$}
$$ N := \{Z\in\C^{n,1} \:  \left\langle Z,Z \right\rangle <0 \} .$$
The \textit{complex hyperbolic n-space $\chs$} is defined to be $\proj(N)$, i.e. the subset of $\proj( \Cnh )$ consisting of negative lines in $\C^{n,1}$. We identify $\chs$ with the open unit ball 
$$ \B^n := \{ z\in\C^n \: (z,z) < 1 \} $$
as follows. Define a mapping $A$ by 
\begin{align*}
     A &: \C^n \longrightarrow \proj(\C^{n,1}) \\
       & z'\longmapsto \left[ \begin{array}{l} z' \\	1	\end{array} \right] \,.
\end{align*}
Since for negative vectors of $\C^{n,1}$ the $(n+1)$-coordinate is necessarily different from zero, $\chs \subset A(\C^n)$. The mapping $A$ identifies $\B^n$ with $\chs$ and $\partial \B^n = \s^{2n-1} \subset \C^n $ with $\partial \chs$.
\begin{theorem}\textnormal{(\cite{goldman99}, Theorem 3.1.10)}
Let $F\subset \proj (\C^{n,1})$ be a complex $m$-di\-men\-sional projective subspace which intersects $\chs$. Then $F \cap \chs$ is a totally geodesic holomorphic submanifold biholomorphically\footnote{biholomorphic mapping = conformal mapping} isometric to $\H_\C^m$.
\end{theorem}
\noindent
The intersection of $\chs$ with a complex hyperplane is a totally geodesic holomorphic complex hypersurface, called a \textit{complex hyperplane in $\chs$}. Its boundary is a smoothly embedded $(2n-3)$-sphere in $\partial \chs$. 

\subsubsection{The Bergman Metric and the Volume Element}\label{subsubsectionTheBergmanMetric}
We normalize the \textit{Berg\-man metric}, a Hermitian metric on $\chs$, to have constant holomorphic sectional curvature $-1$. It can be described as follows. Let $x,y$ be a pair of distinct points in $\B^n$ and let $\overleftrightarrow{xy}$ denote the unique complex line they span. The Bergman metric restricts on $\overleftrightarrow{xy}\cap\B^n$ to the Poincar{\'e} metric of constant curvature $-1$ given by: 
$$ \frac{4 R^2 dz d\overline{z}}{(R^2-r^2)^2} \, ,$$
where $R$ is the radius of the disc $\overleftrightarrow{xy}\cap\B^n$ and $r=r(z)$ is the distance to the center of the disc $\overleftrightarrow{xy}\cap\B^n$. As $\overleftrightarrow{xy}$ is totally geodesic, the distance between $x$ and $y$ in $\chs$ equals the distance between $x$ and $y$ in $\overleftrightarrow{xy}\cap\B^n$ with respect to the above Poincar{\'e} metric. Moreover, the geodesic from $x$ to $y$ in $\chs$ is the Poincar{\'e} geodesic in $\overleftrightarrow{xy}\cap\B^n$ joining $x$ and $y$. The Poincar{\'e} geodesics are circular arcs orthogonal to the boundary and straight lines through the center. 

The volume element on $\chs$ is 
$$d\mu_n = 8^n \frac{r^{2n-1} dr d\sigma}{(1-r^2)^{n+1}} $$
where $d\sigma$ is the volume element on the unit sphere $\s^{2n-1}$. 

\subsection{Convex Geometry}

\subsubsection{Basic Definitions}
The main tool used in this paper is the Fourier transform of distributions, see \cite{GelfandShilov1964} as the classical reference for this topic. As usual, denote by $\sw(\Rn)$ the \textit{Schwartz space} of rapidly decreasing infinitely differentiable functions on $\Rn$, also referred to as \textit{test functions}, and by $\sw'(\Rn)$ the space of \textit{distributions} on $\Rn$, the continuous dual of $\sw(\Rn)$. The Fourier transform $\hat{f}$ of a distribution $f$ is defined by $\left\langle \hat{f}, \varphi \right\rangle = \left\langle f, \hat{\varphi} \right\rangle$ for every test function $\varphi$. A distribution $f$ on $\Rn$ is \textit{even homogeneous of degree $p\in\R$}, if 
$$ \left\langle f(x), \varphi\left(\frac{x}{\alpha}\right) \right\rangle = |\alpha|^{n+p} \left\langle f, \varphi  \right\rangle $$
for every test function $\varphi$ and every $\alpha\in\R, \alpha\neq 0$. The Fourier transform of an even homogeneous distribution of degree $p$ is an even homogeneous distribution of degree $-n-p$. We call a distribution $f$ \textit{positive definite}, if for every test function $\varphi$
$$ \left\langle f(x), \varphi \ast \overline{\varphi}(-x) \right\rangle \geq 0 \, . $$  
This is equivalent to $\hat{f}$ being a positive distribution, i.e. $\left\langle \hat{f}, \varphi \right\rangle \geq 0$ for every non-negative test function $\varphi$.

A compact set $K$ in $\R^n$ is called a \textit{star body} if every line through the origin crosses the boundary in exactly two points different from the origin, and its \textit{Minkowski functional} is defined by
$$ \|x\|_K := \min \{ a\geq 0 \: x\in aK\} \, .$$  
The boundary of $K$ is continuous if $\|\cdot \|_K$ is a continuous function on $\Rn$. If in addition $K$ is origin symmetric and convex, then the Minkowski functional is a norm on $\R^n$. A star body $K$ is said to be \textit{$k$-smooth}, $k\in\N\cup \{0\}$, if the restriction of $\|\cdot \|_K$ to the unit sphere $\sn$ belongs to the class $C^k(\sn)$ of $k$ times continuously differentiable functions on the unit sphere. If $\|\cdot \|_K \in C^k(\sn)$ for any $k\in\N$, then a star body $K$ is said to be \textit{infinitely smooth}. For $x\in\sn$, the \textit{radial function of $K$},  $\rho_K(x)=\|x\|_K^{-1}$, is the Euclidean distance from the origin to the boundary of $K$ in the direction $x$. For all bodies considered in the sequel, the origin is an interior point. 

\subsubsection{Approximation Results}\label{subsubsection_AR}
One can approximate any convex body $K$ in $\Rn$ in the \textit{radial metric} 
$$ \rho(K,L):= \max\limits_{x\in\sn} |\rho_K(x)-\rho_L(x)| $$
by a sequence of infinitely smooth convex bodies with the same symmetries as $K$, see Theorem 3.3.1 in \cite{schneider1993}. In particular, any $R_{\theta}$-invariant convex body in $\R^{2n}$ can be approximated by infinitely smooth $R_{\theta}$-invariant convex bodies. Any $k$-smooth star body $K$ can be approximated by a sequence of infinitely smooth star bodies $K_m$ so that the radial functions $\rho_{K_m}$ converge to $\rho_K$ in the metric of the space $C^k(\sn)$, see \cite{Koldobsky2005}, p. 27, preserving the symmetries of $K$ as well.

\subsubsection{Fourier Approach to Sections}
It was shown in \cite{Koldobsky2005}, Lemma 3.16, that for an infinitely smooth origin sym\-me\-tric star body $K$ in $\Rn$ and $0<p<n$, the Fourier transform of the distribution $\|x\|^{-p}_K$ is an infinitely smooth function on $\Rn \setminus \{0\}$, homogeneous of degree $-n+p$.
We shall use the following analogue of the Parseval's formula:
\begin{lemma}\textnormal{(\cite{Koldobsky2005}, Lemma 3.22)}
Let $K$ and $L$ be infinitely smooth origin symmetric star bodies in $\R^n$, and let $0<p<n$. Then 
$$ \int_{\sn} (\|\cdot\|_K^{-p})^{\wedge}(\theta) (\|\cdot\|_L^{-n+p})^{\wedge}(\theta) d\theta = (2\pi)^n \int_{\sn} \|\theta\|_K^{-p} \|\theta\|_L^{-n+p} d\theta \, . $$
\end{lemma} 
\noindent
We will also use the following version of the above lemma.
\begin{co}\label{coro_PWM}\textnormal{(\cite{Koldobsky2005}, Corollary 3.23)}
Let $k\in \N$ with $0<k<n$. Let $f$ and $g$ be two even functions on $\Rn$, homogeneous of degree $-k$ and $-n+k$, respectively. Suppose that $f$ represents a positive definite distribution on $\Rn$ and let $\mu_0$ be the finite Borel measure on $\sn$ that corresponds to $\hat{f}$ by Corollary 2.26 in \cite{Koldobsky2005}. Then
$$ \int_{\sn} \hat{g}(\theta) \, d\mu_0 = (2\pi)^n \int_{\sn} g(\theta) f(\theta) \, d\theta \, .$$
\end{co} 
Let $0<k<n$ and let $H$ be an $(n-k)$-dimensional subspace of $\Rn$. Fix an orthonormal basis $e_1, \dots, e_k$ in the orthogonal subspace $H^{\perp}$. For a star body $K$ in $\Rn$, define the \textit{$(n-k)$-dimensional parallel section function\footnote{The parallel section function in this generality was first introduced in \cite{Koldobsky2000}, where the author gives a characterization of several classes of generalized intersection bodies in the language of functional analysis. The parallel section function and its further generalizations proved useful for the solution of virtually every generalization of the Busemann-Petty problem, e.g. \cite{KoldobskyYaskinYaskina2006, Zymonopoulou2008}. In particular, the parallel section function on $\R^1$ was used in the Fourier analytic proof of the original Busemann-Petty problem, \cite{GardnerKoldobskySchlumprecht1999}.} $A_{K,H}$} as a function on $\R^k$ such that for $u\in \R^k$
\begin{align*}
	A_{K,H}(u)	&= \Vol_{n-k}(K\cap\{H+u_1 e_1 + \cdots + u_k e_k\}) \\
       				&= \int_{\{ x\in \Rn\: (x,e_1)=u_1, \dots , (x,e_k)=u_k \}} \chi (\|x\|_K) dx\, ,
\end{align*}
where $\chi$ is the indicator function of the interval $[0,1]$. If $K$ is infinitely smooth, the function $A_{K,H}$ is infinitely differentiable at the origin. We shall make use of the following fact:
\begin{lemma}\label{lem_PSFE}\textnormal{(\cite{Koldobsky2000}, Theorem 2)}
Let $K$ be an infinitely smooth origin symmetric star body in $\Rn$ and $0<k<n$. Then for every $(n-k)$-dimensional subspace $H$ of $\Rn$ and for every $m\in \N\cup\{0\}$, $m<(n-k)/2$,
$$ \Delta^m A_{K,H}(0) = \frac{(-1)^m}{(2\pi)^k(n-2m-k)} \int_{\sn\cap H^{\perp}} (\|x\|_K^{-n+2m+k})^{\wedge}(\xi) d\xi \, ,$$
where $\Delta$ denotes the Laplacian on $\R^k$.
\end{lemma}
\subsubsection{The Fourier and Radon Transforms of $R_{\theta}$-invariant Functions}
The following simple observation is crucial for the application of Fourier methods to sections of convex bodies in the complex case. It translates the $R_{\theta}$-invariance of a body $K$ into a certain invariance of the Fourier transform of its Minkowski functional raised to some power. We reproduce the proof of this observation here for its simple and insightful nature.
\begin{lemma}\label{lem_ConD}\textnormal{(\cite{KoldobskyKonigZymonopoulou2008}, Lemma 3)}
Suppose that $K$ is an infinitely smooth $R_{\theta}$-in\-va\-ri\-ant star body in $\R^{2n}$. Then for every $0<p<2n$ and $\xi\in\s^{2n-1}$ the Fourier transform of the distribution $\|x\|_K^{-p}$ is a constant function on $\s^{2n-1}\cap H_{\xi}^{\perp}$.
\end{lemma}
\pf
As mentioned above the Fourier transform of $\|x\|_K^{-p}$ is an in\-fi\-nitely smooth function on $\Rn \setminus \{0\}$. Since the function $\|x\|_K$ is $R_{\theta}$-invariant, by the connection between the Fourier transform of distributions and linear transformations, the Fourier transform of $\|x\|_K^{-p}$ is also $R_{\theta}$-invariant. Recall, from the introduction, that the two-dimensional space $H_{\xi}^{\perp}$ is spanned by two vectors $\xi$ and $\xi_{\perp}$. Every vector in $\s^{2n-1}\cap H_{\xi}^{\perp}$ is the image of $\xi$ under one of the coordinate-wise rotations $R_{\theta}$, so the Fourier transform of $\|x\|_K^{-p}$ is a constant function on $\s^{2n-1}\cap H_{\xi}^{\perp}$. \endpf
Denote by $C_{\theta}(\s^{2n-1})$ the space of $R_{\theta}$-invariant continuous functions on the unit sphere $\s^{2n-1}$, i.e. continuous real-valued functions $f$ satisfying $f(\xi)=f(R_{\theta} \xi)$ for any $\xi\in\s^{2n-1}$ and any $\theta\in [0, 2\pi]$. The \textit{complex spherical Radon transform}, introduced in \cite{KoldobskyPaourisZymonopoulou2011}, is an operator $\RT_c\:C_{\theta}(\s^{2n-1}) \rightarrow C_{\theta}(\s^{2n-1})$ defined by
$$ \RT_c f(\xi) = \int_{\s^{2n-1}\cap H_{\xi}} f(x) dx. $$
To derive a formula for the volume of the $(2n-2)$-dimensional section of an $R_{\theta}$-invariant star body in $\R^{2n}$, contained in the unit ball, by a hyperplane $H_{\xi}$ with respect to the volume element in $\chs$, we use the following recently established connection between the complex spherical Radon transform and the Fourier transform.
\begin{lemma}\textnormal{(\cite{KoldobskyPaourisZymonopoulou2011}, Lemma 4)}
Let $f \in C_{\theta}(\s^{2n-1})$ be an even function. Extend $f$ to a homogeneous function of degree $-2n+2$, $f \cdot r^{-2n+2}$, then the Fourier transform of this extension is a homogeneous function of degree $-2$ on $\R^{2n}$, whose restriction to the unit sphere is continuous. Moreover, for every $\xi\in\s^{2n-1}$
$$ \RT_c f(\xi) = \frac{1}{2\pi} (f \cdot r^{-2n+2})^{\wedge} (\xi) \, . $$
\end{lemma}

\begin{lemma}\label{lem_VS}
Let $K$ be a continuous $R_{\theta}$-invariant star body in $\R^{2n}$ with $n\geq 2$ contained in the unit ball. For $\xi\in\s^{2n-1}$, we have 
\begin{equation}\nonumber
 \HVol_{2n-2}(K\cap H_{\xi}) = \frac{ 8^{n-1}}{2\pi} \left( |x|^{-2n+2} \int_0^{\frac{|x|}{\|x\|_K}} \frac{r^{2n-3}}{(1-r^2)^n} dr \right)^{\wedge} (\xi) \, .
\end{equation}
\end{lemma}

\pf  
We compute:
\begin{align*}
	\HVol_{2n-2}(K\cap H_{\xi}) &= \int_{K\cap H_{\xi}} d\mu_{n-1} \\
       												&= 8^{n-1} \int_{\s^{2n-1}\cap H_{\xi}} \int_0^{\|x\|_K^{-1}} \frac{r^{2n-3}}{(1-r^2)^n} dr dx \, .
\end{align*}
Using that $|x|=1$, we rewrite the above integral as:
$$ \HVol_{2n-2}(K\cap H_{\xi})	= 8^{n-1} \int_{\s^{2n-1}\cap H_{\xi}} |x|^{-2n+2} \int_0^{\frac{|x|}{\|x\|_K}} \frac{r^{2n-3}}{(1-r^2)^n} dr dx \, .$$
The function under the first integral sign is a homogeneous function of degree $-2n+2$ and thus by the above lemma we obtain:
\begin{equation}\nonumber
 \HVol_{2n-2}(K\cap H_{\xi})	= \frac{8^{n-1}}{2\pi} \left( |x|^{-2n+2} \int_0^{\frac{|x|}{\|x\|_K}} \frac{r^{2n-3}}{(1-r^2)^n} dr \right)^{\wedge} (\xi) \, .
\end{equation} \endpf

\section{Connection with the distribution $\frac{\|x\|_K^{{-2}}}{1-\left(\frac{|x|}{\|x\|_K}\right)^2}$}
\noindent
We now turn to our problem. First we construct a counterexample to the Busemann-Petty problem in $\chs$ for $n\geq 4$. We use the same idea as in \cite{yaskin2006}, namely that any Riemannian space is locally close to being Euclidean.   
\begin{theorem}\label{thm_CEIDFOUR}
There exist $R_{\theta}$-invariant convex bodies $K$ and $L$ in $\R^{2n}$ with $n\geq 4$ contained in the unit ball so that
$$ \HVol_{2n-2}(K\cap H_{\xi}) \leq \HVol_{2n-2}(L\cap H_{\xi})$$
for every $\xi \in \s^{2n-1}$, but 
$$ \HVol_{2n}(K) > \HVol_{2n}(L) \, .$$
\end{theorem}

\pf
Let $K$ and $L$ be $R_{\theta}$-invariant convex bodies in $\R^{2n}$ with $n\geq 4$ that provide a counterexample to the complex Busemann-Petty problem, see \cite{KoldobskyKonigZymonopoulou2008}, i.e.
\begin{equation}\label{equationCBPA}
\Vol_{2n-2}(K\cap H_{\xi}) \leq \Vol_{2n-2}(L\cap H_{\xi})
\end{equation}
for every $\xi \in \s^{2n-1}$, but 
\begin{equation}\label{equationCBPCC}
\Vol_{2n}(K) > \Vol_{2n}(L) \, .
\end{equation}
We can and will assume that bodies $K$ and $L$ are infinitely smooth. Note that a dilation by a positive factor is an automorphism of $R_{\theta}$-invariant convex bodies in $\R^{2n}$. Since inequality (\ref{equationCBPCC}) is strict, we can dilate the body $L$ by a positive factor greater than one, to make inequality (\ref{equationCBPA}) strict as well. Furthermore, there is an $\epsilon > 0$ so that 
$$ (1+\epsilon) \Vol_{2n-2}(K\cap H_{\xi}) \leq \Vol_{2n-2}(L\cap H_{\xi})$$
for every $\xi \in \s^{2n-1}$, but 
$$  \Vol_{2n}(K) > (1+\epsilon) \Vol_{2n}(L) \, .$$
Here we are using the fact that $\xi \mapsto A_{K,H_{\xi}}(0)$ is continuous. It follows by an argument similar to one given in \cite{Koldobsky2005}, Lemma 2.4. And any dilation of bodies $K$ and $L$ by a factor $\alpha>0$ will also provide a counterexample. 

Choose $\alpha$ so small that both bodies $\alpha K$ and $\alpha L$ lie in the ball of radius $s$ that satisfies the inequality
$$ 1 \leq \frac{1}{(1-s^2)^{n+1}} \leq 1+\epsilon \, .$$
We use the same letters for the dilated bodies. Then for the volumes of bodies $K$ and $L$ we obtain:
\begin{align*}
	\HVol_{2n}(L) 
       					&= 8^n \int_L \frac{dx_{2n}}{(1-|x|^2)^{n+1}} \\
       					&\leq 8^n (1+\epsilon) \int_L dx_{2n} \\
       					&= 8^n (1+\epsilon) \Vol_{2n}(L) \\
       					&< 8^n \Vol_{2n}(K) \\
       					&= 8^n \int_K dx_{2n} \\
       					&\leq 8^n \int_K \frac{dx_{2n}}{(1-|x|^2)^{n+1}}  \\
       					&= \HVol_{2n}(K) \, ,
\end{align*}
and similarly, for the volumes of sections we have
\begin{align*}
	\HVol_{2n-2}(K\cap H_{\xi})	&= 8^{n-1} \int_{K\cap H_{\xi}} \frac{dx_{2n-2}}{(1-|x|^2)^n} \\
       												&\leq 8^{n-1} (1+\epsilon) \int_{K\cap H_{\xi}} dx_{2n-2} \\
       												&= 8^{n-1} (1+\epsilon) \Vol_{2n-2}({K\cap H_{\xi}}) \\
       												&\leq 8^{n-1} \Vol_{2n-2}({L\cap H_{\xi}}) \\
       												&\leq 8^{n-1} \int_{L\cap H_{\xi}} \frac{dx_{2n-2}}{(1-|x|^2)^{n}}  \\
       												&= \HVol_{2n-2}({L\cap H_{\xi}}) \, .
\end{align*} \endpf

The connection between the Busemann-Petty problem in $\chs$ and the distribution $\frac{\|x\|_K^{{-2}}}{1-\left(\frac{|x|}{\|x\|_K}\right)^2}$ is the following:
\begin{theorem}\label{thm_CbeAtoBRPandD}
The answer to the Busemann-Petty problem in $\chs$ is affirmative if and only if for every $R_{\theta}$-invariant convex body in $\R^{2n}$ contained in the unit ball the distribution $\frac{\|x\|_K^{{-2}}}{1-\left(\frac{|x|}{\|x\|_K}\right)^2}$ is positive definite.
\end{theorem}
This theorem will follow from the next two lemmas.
\begin{lemma}
Let $K$ be an $R_{\theta}$-invariant convex body in $\R^{2n}$ contained in the unit ball such that $\frac{\|x\|_K^{{-2}}}{1-\left(\frac{|x|}{\|x\|_K}\right)^2}$ is a positive definite distribution on $\R^{2n}$.
And let $L$ be an $R_{\theta}$-invariant star body in $\R^{2n}$ contained in the unit ball so that 
$$ \HVol_{2n-2}(K\cap H_{\xi}) \leq \HVol_{2n-2}(L\cap H_{\xi})$$
for every $\xi \in \s^{2n-1}$. Then  
$$ \HVol_{2n}(K) \leq \HVol_{2n}(L) \, .$$
\end{lemma}
\pf
Observe that function $\frac{r^2}{1-r^2}$ is an increasing function on the interval $(0,1)$. We use this observation to estimate the following expression:
\begin{align*}
	\frac{a^2}{1-a^2}\int_a^b \frac{r^{2n-3}}{(1-r^2)^n} dr 
						&= \frac{a^2}{1-a^2}\int_a^b \frac{r^{2n-1}}{(1-r^2)^{n+1}} \frac{(1-r^2)}{r^2} dr \\
       			&= \int_a^b \frac{r^{2n-1}}{(1-r^2)^{n+1}} \frac{a^2}{1-a^2} \frac{(1-r^2)}{r^2} dr \\
       			&\leq \int_a^b \frac{r^{2n-1}}{(1-r^2)^{n+1}} dr \, ,
\end{align*}
where $a, b$ are in $(0,1)$. Observe that the above inequality it true in case $a\leq b$ as well as in case $b\leq a$.
Integrating both sides in the above inequality over the unit sphere $\s^{2n-1}$ with $a=\|x\|^{-1}_K$ and $b=\|x\|^{-1}_L$ we obtain:
\begin{equation}\label{eq_AoS}
\int_{\s^{2n-1}} 	\frac{\|x\|^{-2}_K}{1-\|x\|^{-2}_K} \int_{\|x\|^{-1}_K}^{\|x\|^{-1}_L} \frac{r^{2n-3}}{(1-r^2)^n} dr dx 	\leq \int_{\s^{2n-1}} \int_{\|x\|^{-1}_K}^{\|x\|^{-1}_L} \frac{r^{2n-1}}{(1-r^2)^{n+1}} dr dx \, .
\end{equation}
Next we show that the left hand side in the above expression is positive. This amounts to showing that 
$$ \int_{\s^{2n-1}} 	\frac{\|x\|^{-2}_K}{1-\|x\|^{-2}_K} \int^{\|x\|^{-1}_K}_0 \frac{r^{2n-3}}{(1-r^2)^n} dr dx 	\leq \int_{\s^{2n-1}} 	\frac{\|x\|^{-2}_K}{1-\|x\|^{-2}_K}  \int_0^{\|x\|^{-1}_L} \frac{r^{2n-3}}{(1-r^2)^n} dr dx \, . $$
Indeed, let $d\mu_0$ be the measure corresponding to the Fourier transform of the positive definite distribution $\frac{\|x\|_K^{{-2}}}{1-\left(\frac{|x|}{\|x\|_K}\right)^2}$, then using Corollary \ref{coro_PWM} we obtain:
\begin{align*}
	(2\pi)^{2n} & \int_{\s^{2n-1}} 	\frac{\|x\|^{-2}_K}{1-\|x\|^{-2}_K} \int^{\|x\|^{-1}_K}_0 \frac{r^{2n-3}}{(1-r^2)^n} dr dx \\
						&= (2\pi)^{2n} \int_{\s^{2n-1}} \left(	\frac{\|x\|^{-2}_K}{1-\left( \frac{|x|}{\|x\|_K} \right)^2 } \right) \left( |x|^{-2n+2} \int^{\frac{|x|}{\|x\|_K}}_0 \frac{r^{2n-3}}{(1-r^2)^n} dr \right) dx \\
       			&= \int_{\s^{2n-1}} \left( |x|^{-2n+2} \int^{\frac{|x|}{\|x\|_K}}_0 \frac{r^{2n-3}}{(1-r^2)^n} dr \right)^{\wedge} (\xi) \, d\mu_0(\xi) \\
       			&= \frac{2 \pi}{8^{n-1}} \int_{\s^{2n-1}} \HVol_{2n-2}(K\cap H_{\xi}) \, d\mu_0(\xi) \\
       			&\leq \frac{2 \pi}{8^{n-1}} \int_{\s^{2n-1}} \HVol_{2n-2}(L\cap H_{\xi}) \, d\mu_0(\xi) \\
       			&= \int_{\s^{2n-1}} \left( |x|^{-2n+2} \int^{\frac{|x|}{\|x\|_L}}_0 \frac{r^{2n-3}}{(1-r^2)^n} dr \right)^{\wedge} (\xi) \, d\mu_0(\xi) \\
       			&= (2\pi)^{2n} \int_{\s^{2n-1}} 	\frac{\|x\|^{-2}_K}{1-\|x\|^{-2}_K} \int^{\|x\|^{-1}_L}_0 \frac{r^{2n-3}}{(1-r^2)^n} dr dx  \, ,
\end{align*}
This implies that the right hand side in equation (\ref{eq_AoS}) is positive as well, which, in turn, shows 
$$ \int_{\s^{2n-1}} \int^{\|x\|^{-1}_K}_0 \frac{r^{2n-1}}{(1-r^2)^{n+1}} dr dx 	\leq \int_{\s^{2n-1}} \int_0^{\|x\|^{-1}_L} \frac{r^{2n-1}}{(1-r^2)^{n+1}} dr dx \, . $$
That is, 
$$ \HVol_{2n}(K) \leq \HVol_{2n}(L) \, . \vspace{-3.29mm}$$ \endpf
Next result is the inversion of the previous lemma. Its proof is based on a standard perturbation argument, see \cite{Koldobsky1999}, Theorem 2, \cite{Zvavitch2005}, Theorem 2, or \cite{yaskin2006}, Proposition 3.9.
\begin{lemma}\label{lem_CE}
Suppose there is an infinitely smooth $R_{\theta}$-invariant convex body $K$ contained in the unit ball of $\R^{2n}$ with strictly positive curvature\footnote{By strictly positive curvature we mean that the boundary of $K$ does not contain any straight line segments. More precisely, the normal curvature of $K$ at any point $p$ on the boundary of $K$ is strictly positive in any direction $v$, where $v\in TK_p$ is any element in the tangent space to $K$ at $p$, see \cite{Thorpe79}.} so that $\frac{\|x\|_K^{{-2}}}{1-\left(\frac{|x|}{\|x\|_K}\right)^2}$ is not a positive definite distribution on $\R^{2n}$. Then one can perturb the body $K$ to construct another $R_{\theta}$-invariant convex body $L$ contained in $\s^{2n-1}$ so that for every $\xi \in \s^{2n-1}$ 
$$ \HVol_{2n-2}(L\cap H_{\xi}) \leq \HVol_{2n-2}(K\cap H_{\xi}) \, ,$$
but
$$ \HVol_{2n}(L) > \HVol_{2n}(K) \, .$$  
\end{lemma}

\pf
It follows from our assumptions that the Fourier transform of $\frac{\|x\|_K^{{-2}}}{1-\left(\frac{|x|}{\|x\|_K}\right)^2}$ is negative on some open subset $\Omega$ of the sphere $\s^{2n-1}$. $R_{\theta}$-invariance of the body $K$ implies the $R_{\theta}$-invariance of the set $\Omega$. Choose a smooth non-negative $R_{\theta}$-invariant function $f$ on $\s^{2n-1}$ with support of $f$ contained in $\Omega$ and extend $f$ to an $R_{\theta}$-invariant homogeneous function $f\left(\frac{x}{|x|}\right) |x|^{-2}$ of degree $-2$ on $\R^{2n}$. The Fourier transform of this extension is an $R_{\theta}$-invariant infinitely smooth function on $\R^{2n}\setminus\{0\}$,  homogeneous of degree $-2n+2$, i.e.
$$ \left( f\left(\frac{x}{|x|}\right) |x|^{-2} \right)^{\wedge} (y) = g\left(\frac{y}{|y|}\right) |y|^{-2n+2} $$
with $g\in C^{\infty}(\s^{2n-1})$. Since $f$ is $R_{\theta}$-invariant, so is $g$. Define an origin symmetric body $L$ contained in $\s^{2n-1}$ by
$$ |x|^{-2n+2} \int^{\frac{|x|}{\|x\|_L}}_0 \frac{r^{2n-3}}{(1-r^2)^{n}} dr = |x|^{-2n+2} \int^{\frac{|x|}{\|x\|_K}}_0 \frac{r^{2n-3}}{(1-r^2)^{n}} dr -\epsilon g\left(\frac{x}{|x|}\right) |x|^{-2n+2} $$
for some $\epsilon >0$. By a similar argument as in \cite{Zvavitch2005}, Proposition 2, for small enough $\epsilon$, the body $L$ is convex. From the $R_{\theta}$-invariance of $K$ and $g$ follows the $R_{\theta}$-invariance of $L$. Thus $L$ is an $R_{\theta}$-invariant convex body contained in $\s^{2n-1}$. Using Lemma \ref{lem_VS}, we compute: 
\begin{align*}
	\HVol_{2n-2}(L\cap H_{\xi}) &= \frac{8^{n-1}}{2 \pi} \left( |x|^{-2n+2} \int^{\frac{|x|}{\|x\|_L}}_0 \frac{r^{2n-3}}{(1-r^2)^n} dr \right)^{\wedge} (\xi) \\
															&= \frac{ 8^{n-1}}{2 \pi} \left( |x|^{-2n+2} \int^{\frac{|x|}{\|x\|_K}}_0 \frac{r^{2n-3}}{(1-r^2)^n} dr \right)^{\wedge} (\xi) \\
															& \,\,\,\,\,\,\,\,\,\,\,\,\,\,\,\,\,\,\,\,\,\,\,\,\,\,\,\,\,\,\,\,\,\,\,\,\,\,\,\,\,\,\,\,\,\,\,\,\,\,\,\,\,\,\,\,\,\,\,\,\,\,\,\,\,\,\,\,\,\,\,\,\,\,\,\,\,\,\,\,\,\,\,\,\,\,\,\,\,\,\,\,\, - \epsilon (2\pi)^{2n} f\left(\frac{x}{|x|}\right) |x|^{-2} \\
															&\leq \frac{ 8^{n-1}}{2 \pi} \left( |x|^{-2n+2} \int^{\frac{|x|}{\|x\|_K}}_0 \frac{r^{2n-3}}{(1-r^2)^n} dr \right)^{\wedge} (\xi) \\
       												&= \HVol_{2n-2}(K\cap H_{\xi}).
\end{align*}

As in the proof of the previous lemma, to complete the proof it is enough to show the following: 
\begin{align*}
	(2\pi)^{2n} & \int_{\s^{2n-1}} 	\frac{\|x\|^{-2}_K}{1-\|x\|^{-2}_K} \int^{\|x\|^{-1}_L}_0 \frac{r^{2n-3}}{(1-r^2)^n} dr dx \\
							&= (2\pi)^{2n} \int_{\s^{2n-1}} \left(	\frac{\|x\|^{-2}_K}{1-\left( \frac{|x|}{\|x\|_K} \right)^2 } \right) \left( |x|^{-2n+2} \int^{\frac{|x|}{\|x\|_L}}_0 \frac{r^{2n-3}}{(1-r^2)^n} dr \right) dx \\
							&= \int_{\s^{2n-1}} \left( 	\frac{\|x\|^{-2}_K}{1-\left( \frac{|x|}{\|x\|_L} \right)^2 } \right)^{\wedge} (\xi) \left( |x|^{-2n+2} \int^{\frac{|x|}{\|x\|_L}}_0 \frac{r^{2n-3}}{(1-r^2)^n} dr \right)^{\wedge} (\xi) \, d\xi \\
       				&= \int_{\s^{2n-1}} \left( 	\frac{\|x\|^{-2}_K}{1-\left( \frac{|x|}{\|x\|_K} \right)^2 } \right)^{\wedge} (\xi) \left( |x|^{-2n+2} \int^{\frac{|x|}{\|x\|_K}}_0 \frac{r^{2n-3}}{(1-r^2)^n} dr \right)^{\wedge} (\xi) \, d\xi \\
       				& \,\,\,\,\,\,\,\,\,\,\,\,\,\,\,\,\,\,\,\,\,\,\,\,\,\,\,\,\,\,\,\,\,\,\,\,\,\,\,\,\,\,\,\,\,\,\,\,\,\,\,\,\,\,\,\,\,\,\,\,\,\, -\epsilon (2\pi)^{2n} \int_{\s^{2n-1}}  \left( 	\frac{\|x\|^{-2}_K}{1-\left( \frac{|x|}{\|x\|_K} \right)^2 } \right)^{\wedge} (\xi) f(\xi) d\xi \\				
       				&> \int_{\s^{2n-1}} \left( 	\frac{\|x\|^{-2}_K}{1-\left( \frac{|x|}{\|x\|_K} \right)^2 } \right)^{\wedge} (\xi) \left( |x|^{-2n+2} \int^{\frac{|x|}{\|x\|_K}}_0 \frac{r^{2n-3}}{(1-r^2)^n} dr \right)^{\wedge} (\xi) \, d\xi \\
       				&= (2\pi)^{2n} \int_{\s^{2n-1}} 	\frac{\|x\|^{-2}_K}{1-\|x\|^{-2}_K} \int^{\|x\|^{-1}_K}_0 \frac{r^{2n-3}}{(1-r^2)^n} dr dx \, . 
\end{align*} \endpf

\section{Solution of the Problem}
\noindent
In view of Theorem \ref{thm_CEIDFOUR} we only have to find out whether the distribution $\frac{\|x\|_K^{{-2}}}{1-\left(\frac{|x|}{\|x\|_K}\right)^2}$ is positive definite for the complex dimension two and three, the case of the complex dimension one is trivial.\footnote{Indeed, $R_{\theta}$-invariant convex bodies in $\C=\R^2$ are closed discs around the origin.} The parallel section function, defined earlier, will help us to carry out this task. The idea is to express the distribution of interest in terms of the parallel section function.

In dimension $2n$, Lemma \ref{lem_PSFE} with $H=H_{\xi}$, $\xi \in \s^{2n-1}$, in which case $k=2$, reads as follows: Let $K$ be an infinitely smooth $R_{\theta}$-invariant star body in $\R^{2n}$, then for $m\in \N\cup \{0\}$, $m<n-1$
$$ \Delta^m A_{K,H_{\xi}}(0) = \frac{(-1)^m}{(2\pi)^2 (2n-2m-2)} \int_{\s^{2n-1}\cap H_{\xi}^{\perp}} (\|x\|_K^{-2n+2m+2})^{\wedge}(\nu) d\nu \, .$$
Since the above integral is taken over the region $\s^{2n-1}\cap H_{\xi}^{\perp}$, by Lemma \ref{lem_ConD}, it follows that
\begin{equation}\label{eq_PSF}
		\Delta^m A_{K,H_{\xi}}(0)=\frac{(-1)^m}{2\pi(2n-2m-2)} (\|x\|_K^{-2n+2m+2})^{\wedge}(\xi) \, .
\end{equation}
In complex dimension two there is only one choice for $m$, namely $m=0$, and equation (\ref{eq_PSF}) becomes
$$ A_{K,H_{\xi}}(0) = \frac{1}{4\pi} (\|x\|_K^{-2})^{\wedge}(\xi)\, . $$
In complex dimension three, evaluating equation (\ref{eq_PSF}) for $m=1$, we obtain
$$ 	\Delta A_{K,H_{\xi}}(0) = \frac{-1}{4\pi} (\|x\|_K^{-2})^{\wedge}(\xi) \, . $$

\begin{lemma}\label{lem_DPD}
For any $R_{\theta}$-invariant star body $K$ in $\R^4$ contained in the unit ball the distribution $\frac{\|x\|_K^{{-2}}}{1-\left(\frac{|x|}{\|x\|_K}\right)^2}$ is positive definite.
\end{lemma}
\pf
Assume first that $K$ is an infinitely smooth $R_{\theta}$-invariant star body in $\R^4$. Define another body $M$ by
$$ \|x\|^{-2}_M = \frac{\|x\|_K^{{-2}}}{1-\left(\frac{|x|}{\|x\|_K}\right)^2} \, . $$   
Since $K$ is an $R_{\theta}$-invariant star body, so is $M$. Note also that the denominator in the defining expression for $M$ is never zero, since the body $K$ is contained in the open unit ball. This implies that the body $M$ is also infinitely smooth. Hence
$$ A_{M,H_{\xi}}(0) = \frac{1}{4\pi} (\|x\|_M^{-2})^{\wedge}(\xi)\, , $$
and consequently $\|x\|_M^{-2}$ is positive definite. Thus for any infinitely smooth $R_{\theta}$-invariant star body $K$ the distribution $\frac{\|x\|_K^{{-2}}}{1-\left(\frac{|x|}{\|x\|_K}\right)^2}$ is positive definite. And therefore, by the approximation results in Section \ref{subsubsection_AR}, this is true for any $R_{\theta}$-invariant star body. \endpf
\begin{lemma}\label{lem_CEFNthree}
There is an infinitely smooth $R_{\theta}$-invariant convex body $K$ in $\R^6$ contained in the unit ball for which the distribution $\frac{\|x\|_K^{{-2}}}{1-\left(\frac{|x|}{\|x\|_K}\right)^2}$ is not positive definite.
\end{lemma}

\pf
For an element $\xi$ of $\R^6$, $\xi=(\xi_{11}, \xi_{12}, \xi_{21}, \xi_{22}, \xi_{31}, \xi_{32})$, denote by $\xi_3=(\xi_{31}, \xi_{32})$ and by $\tilde{\xi}=(\xi_{11}, \xi_{12}, \xi_{21}, \xi_{22})$, then $\xi=(\tilde{\xi}, \xi_3)$. We will work with the following map, written in polar coordinates as 
\begin{equation}\label{eq_MinPC}
(r,\theta) \mapsto \left( \sqrt{\frac{r^2}{1-r^2}}, \theta \right) \, .
\end{equation}
Note that this map, restricted to the two-dimensional plane $xy$, takes the line $x=\frac{1}{a}$ to the hyperbola $(a^2-1) x^2 - y^2 =1$, provided that $a^2-1>0$, and it takes the ellipse $x^2 + (1+b^2) y^2 =1$ to the line $y=\frac{1}{b}$.\footnote{Indeed, writing the equation of the line $x=\frac{1}{a}$ in polar coordinates, we obtain $\frac{r^2}{1-r^2} = \frac{1}{a^2 \cos^2 \theta -1}$, and so the image of the line is $r^2 = \frac{1}{a^2 \cos^2 \theta -1}$. Similarly to find the image of the ellipse, we write its equation in polar coordinates: $r^2 =\frac{1}{1 + b^2 \sin^2 \theta}$. Then $\frac{r^2}{1-r^2} = \frac{1}{b^2 \sin^2 \theta}$ and hence the image of the ellipse is $r=\frac{1}{b \sin \theta}$.} Denote the equation of the elliptic arc above the $x$-axis by $e$, i.e. $e(x)=\sqrt{\frac{1-x^2}{1+b^2}}$, and the equation of the hyperbolic arc to the right of the $y$-axis by $h$, i.e. $h(y)=\sqrt{\frac{1+y^2}{a^2-1}}$. Now define a convex body $K$ in $\R^6$ by
$$ K=\{\xi\in \R^6 \: |\tilde{\xi}| \leq \frac{1}{a} \text{ and } |\xi_3|\leq e(|\tilde{\xi}|)\} \, .$$
We restrict the values of $b$ to be strictly greater then one, this ensures that the body $K$ is contained in the unit ball. As before, define a star body $M$ by
$$ \|x\|^{-2}_M = \frac{\|x\|_K^{{-2}}}{1-\left(\frac{|x|}{\|x\|_K}\right)^2} \, . $$  
The body $M$ is an image of the body $K$ under the map (\ref{eq_MinPC}) and hence it can be described as 
$$ M=\{\xi\in \R^6 \: |\tilde{\xi}| \leq h(|\xi_3|) \text{ with } |\xi_3|\leq \frac{1}{b}\} \, .$$
The bodies $K$ and $M$ we constructed are not infinitely smooth. However, since we can approximate the body $K$ in the radial metric by a sequence of infinitely smooth convex bodies, there is an infinitely smooth convex body $K'$ that differs from $K$ only in an arbitrary small neighborhood of the boundary of $K$. This modification of $K$ will make the body $M$ infinitely smooth as well. We use the same letters $K$ and $M$ for these modified bodies. By construction both bodies are $R_{\theta}$-invariant.

Let $x=(\tilde{x}, x_3)\in M$ with $x_3\neq (0,0)$. Choose $\xi\in \s^5$ in the direction of $x_3$. Fix an orthonormal basis $\{e_1, e_2\}$ for $H_{\xi}^{\perp}$. For $u\in\R^2$ with $|u|>\frac{1}{b}$, $A_{M,H_{\xi}}(u)=0$, and otherwise 
\begin{align*}
	A_{M,H_{\xi}}(u) 	&= \Vol_4 (M\cap\{H_{\xi} + u_1 e_1 + u_2 e_2 \}) \\
										&= \int_{\{ x\in \R^{2n}\: (x,e_1)=u_1, (x,e_2)=u_2 \}} \chi (\|x\|_M) dx \\
										&= \int_{\s^{3}} \int_0^{h(|u|)} r^{3} dr d\theta \\
										&= |\s^{3}| \, \frac{h(|u|)^{4}}{4} \\		
										&= \frac{\pi^2}{2} \, h(|u|)^4 \, ,
\end{align*}
where $|\s^{n-1}|$ stands for the surface area of the unit sphere $\sn$ in $\Rn$, i.e. $2 \pi^{\frac{n}{2}}/\Gamma (\frac{n}{2})$. Setting $a=2$, we get
$$ A_{M,H_{\xi}}(u) = \frac{\pi^2}{2} \left( \frac{1+|u|^2}{3}\right)^2 \text{ and consequently } \Delta A_{M,H_{\xi}}(u) = \frac{4 \pi^2}{9} (1+2|u|^2) \, .$$
Since $M$ is infinitely smooth we have 
$$ (\|x\|_M^{-2})^{\wedge}(\xi)= - 4\pi	\, \Delta A_{M,H_{\xi}}(0) = - \frac{16 \pi^3}{9} \, . $$
This shows that $ \left( \frac{\|x\|_K^{{-2}}}{1-\left(\frac{|x|}{\|x\|_K}\right)^2} \right)^{\wedge}(\xi) = (\|x\|_M^{-2})^{\wedge}(\xi) $ is negative for some direction $\xi$. \endpf

Now we are ready to prove the main result of this paper:
\begin{theorem}
The answer to the Busemann-Petty problem in the complex hyperbolic $n$-space, $\chs$, is affirmative for $n\leq 2$ and negative for $n\geq 3$. 
\end{theorem}

\pf
By Lemma \ref{lem_DPD}, the distribution $\frac{\|x\|_K^{{-2}}}{1-\left(\frac{|x|}{\|x\|_K}\right)^2}$ is positive definite for any $R_{\theta}$-invariant convex body $K$ in $\R^4$, as any such body is a star body. The affirmative answer for $n=2$ now follows from Theorem \ref{thm_CbeAtoBRPandD}.

For $n\geq 4$ the negative answer was provided in Theorem \ref{thm_CEIDFOUR}.
For $n=3$, by Lemma \ref{lem_CEFNthree} there is an infinitely smooth $R_{\theta}$-invariant convex body $K$ in $\R^6$ contained in the unit ball for which the distribution $\frac{\|x\|_K^{{-2}}}{1-\left(\frac{|x|}{\|x\|_K}\right)^2}$ is not positive definite. Observe that we can assume that the body $K$ has strictly positive curvature by applying a standard perturbation trick, namely by setting
$$ \| x \|_{K'}^{-1} = \| x \|_{K}^{-1} + \epsilon |x| $$
for some $\epsilon >0$ small. The negative answer now follows from Lemma \ref{lem_CE}. \endpf


\section*{Acknowledgements}
The author wishes to thank A. Koldobsky for proposing to work on this pro\-blem, for encouragement and many suggestions as well as to V. Yaskin for useful discussions and to the anonymous referee for the thorough review of the paper and many constructive suggestions. 

\nocite{}
\bibliographystyle{amsplain}
\bibliography{ref_01}

\def\dbar{\leavevmode\hbox to 0pt{\hskip.2ex \accent"16\hss}d}
  \def\dbar{\leavevmode\hbox to 0pt{\hskip.2ex \accent"16\hss}d}
  \def\cprime{$'$} \def\dbar{\leavevmode\hbox to 0pt{\hskip.2ex
  \accent"16\hss}d}
\providecommand{\bysame}{\leavevmode\hbox to3em{\hrulefill}\thinspace}
\providecommand{\MR}{\relax\ifhmode\unskip\space\fi MR }
\providecommand{\MRhref}[2]{%
  \href{http://www.ams.org/mathscinet-getitem?mr=#1}{#2}
}
\providecommand{\href}[2]{#2}
\begin{thebibliography}{10}

\bibitem{AbardiaBernig2011}
Judit Abardia and Andreas Bernig, \emph{Projection bodies in complex vector
  spaces}, Adv. Math. \textbf{227} (2011), no.~2, 830--846.

\bibitem{AbardiaGallego2011}
Judit Abardia and Eduardo Gallego, \emph{Convexity on complex hyperbolic
  space}, arXiv:1003. 4667.

\bibitem{Ba}
Keith Ball, \emph{Some remarks on the geometry of convex sets}, Geometric
  aspects of functional analysis (1986/87), Lecture Notes in Math., vol. 1317,
  Springer, Berlin, 1988, pp.~224--231.

\bibitem{BezdekSchneider2010}
K{\'a}roly Bezdek and Rolf Schneider, \emph{Covering large balls with convex
  sets in spherical space}, Beitr\"age Algebra Geom. \textbf{51} (2010), no.~1,
  229--235.

\bibitem{Bo}
Jean Bourgain, \emph{On the {B}usemann-{P}etty problem for perturbations of the
  ball}, Geom. Funct. Anal. \textbf{1} (1991), no.~1, 1--13.

\bibitem{BourgainZhang1999}
Jean Bourgain and Gaoyong Zhang, \emph{On a generalization of the
  {B}usemann-{P}etty problem}, Convex geometric analysis ({B}erkeley, {CA},
  1996), Math. Sci. Res. Inst. Publ., vol.~34, Cambridge Univ. Press,
  Cambridge, 1999, pp.~65--76.

\bibitem{BusemannPetty1956}
Herbert Busemann and Clinton~M. Petty, \emph{Problems on convex bodies}, Math.
  Scand. \textbf{4} (1956), 88--94.

\bibitem{GaoHugSchneider2001}
Fuchang Gao, Daniel Hug, and Rolf Schneider, \emph{Intrinsic volumes and polar
  sets in spherical space}, Math. Notae \textbf{41} (2001/02), 159--176 (2003),
  Homage to Luis Santal{\'o}. Vol. 1 (Spanish).

\bibitem{Ga1}
Richard~J. Gardner, \emph{Intersection bodies and the {B}usemann-{P}etty
  problem}, Trans. Amer. Math. Soc. \textbf{342} (1994), no.~1, 435--445.

\bibitem{Ga2}
\bysame, \emph{A positive answer to the {B}usemann-{P}etty problem in three
  dimensions}, Ann. of Math. (2) \textbf{140} (1994), no.~2, 435--447.

\bibitem{Gardner2002}
\bysame, \emph{The {B}runn-{M}inkowski inequality}, Bull. Amer. Math. Soc.
  (N.S.) \textbf{39} (2002), no.~3, 355--405.

\bibitem{GardnerKoldobskySchlumprecht1999}
Richard~J. Gardner, Alexander Koldobsky, and Thomas Schlumprecht, \emph{An
  analytic solution to the {B}usemann-{P}etty problem on sections of convex
  bodies}, Ann. of Math. (2) \textbf{149} (1999), no.~2, 691--703.

\bibitem{GelfandShilov1964}
Israel~M. Gel{\cprime}fand and Georgi~E. Shilov, \emph{Generalized functions.
  {V}ol. 1}, Academic Press [Harcourt Brace Jovanovich Publishers], New York,
  1964 [1977], Properties and operations, Translated from the Russian by Eugene
  Saletan.

\bibitem{Gi}
Apostolos~A. Giannopoulos, \emph{A note on a problem of {H}. {B}usemann and
  {C}. {M}. {P}etty concerning sections of symmetric convex bodies},
  Mathematika \textbf{37} (1990), no.~2, 239--244.

\bibitem{goldman99}
William~M. Goldman, \emph{Complex hyperbolic geometry}, Oxford Mathematical
  Monographs, The Clarendon Press Oxford University Press, New York, 1999,
  Oxford Science Publications.

\bibitem{K1}
Alexander Koldobsky, \emph{An application of the {F}ourier transform to
  sections of star bodies}, Israel J. Math. \textbf{106} (1998), 157--164.

\bibitem{K2}
\bysame, \emph{Intersection bodies, positive definite distributions, and the
  {B}usemann-{P}etty problem}, Amer. J. Math. \textbf{120} (1998), no.~4,
  827--840.

\bibitem{Koldobsky1999}
\bysame, \emph{A generalization of the {B}usemann-{P}etty problem on sections
  of convex bodies}, Israel J. Math. \textbf{110} (1999), 75--91.

\bibitem{Koldobsky2000}
\bysame, \emph{A functional analytic approach to intersection bodies}, Geom.
  Funct. Anal. \textbf{10} (2000), no.~6, 1507--1526.

\bibitem{Koldobsky2002}
\bysame, \emph{On the derivatives of {X}-ray functions}, Arch. Math. (Basel)
  \textbf{79} (2002), no.~3, 216--222.

\bibitem{Koldobsky2003}
\bysame, \emph{The {B}usemann-{P}etty problem via spherical harmonics}, Adv.
  Math. \textbf{177} (2003), no.~1, 105--114.

\bibitem{Koldobsky2005}
\bysame, \emph{Fourier analysis in convex geometry}, Mathematical Surveys and
  Monographs, vol. 116, American Mathematical Society, Providence, RI, 2005.

\bibitem{koldobsky2011}
\bysame, \emph{Stability of volume comparison for complex convex bodies}, Arch.
  Math. (Basel) \textbf{97} (2011), no.~1, 91--98. \MR{2820591}

\bibitem{KoldobskyKonigZymonopoulou2008}
Alexander Koldobsky, Hermann K{\"o}nig, and Marisa Zymonopoulou, \emph{The
  complex {B}usemann-{P}etty problem on sections of convex bodies}, Adv. Math.
  \textbf{218} (2008), no.~2, 352--367.

\bibitem{KoldobskyPaourisZymonopoulou2011}
Alexander Koldobsky, Grigoris Paouris, and Marisa Zymonopoulou, \emph{Complex
  intersection bodies}, preprint.

\bibitem{KoldobskyYaskinYaskina2006}
Alexander Koldobsky, Vladyslav Yaskin, and Maryna Yaskina, \emph{Modified
  {B}usemann-{P}etty problem on sections of convex bodies}, Israel J. Math.
  \textbf{154} (2006), 191--207.

\bibitem{LR}
David~G. Larman and C.~Ambrose Rogers, Rogers, \emph{The existence of a
  centrally symmetric convex body with central sections that are unexpectedly
  small}, Mathematika \textbf{22} (1975), no.~2, 164--175.

\bibitem{Lu}
Erwin Lutwak, \emph{Intersection bodies and dual mixed volumes}, Adv. in Math.
  \textbf{71} (1988), no.~2, 232--261.

\bibitem{Pa}
Michael Papadimitrakis, \emph{On the {B}usemann-{P}etty problem about convex,
  centrally symmetric bodies in {$\bold R^n$}}, Mathematika \textbf{39} (1992),
  no.~2, 258--266.

\bibitem{rubin2010}
Boris Rubin, \emph{Comparison of volumes of convex bodies in real, complex, and
  quaternionic spaces}, Adv. Math. \textbf{225} (2010), no.~3, 1461--1498.
  \MR{2673737 (2011h:52010)}

\bibitem{schneider1993}
Rolf Schneider, \emph{Convex bodies: the {B}runn-{M}inkowski theory},
  Encyclopedia of Mathematics and its Applications, vol.~44, Cambridge
  University Press, Cambridge, 1993.

\bibitem{Thorpe79}
John~A. Thorpe, \emph{Elementary topics in differential geometry},
  Springer-Verlag, New York, 1979, Undergraduate Texts in Mathematics.

\bibitem{yaskin2006}
Vladyslav Yaskin, \emph{The {B}usemann-{P}etty problem in hyperbolic and
  spherical spaces}, Adv. Math. \textbf{203} (2006), no.~2, 537--553.

\bibitem{Zh1}
Gao~Yong Zhang, \emph{Intersection bodies and the {B}usemann-{P}etty
  inequalities in {${\bf R}^4$}}, Ann. of Math. (2) \textbf{140} (1994), no.~2,
  331--346. \MR{1298716 (95i:52004)}

\bibitem{Zh2}
\bysame, \emph{A positive solution to the {B}usemann-{P}etty problem in {$\bold
  R^4$}}, Ann. of Math. (2) \textbf{149} (1999), no.~2, 535--543.

\bibitem{Zvavitch2005}
Artem Zvavitch, \emph{The {B}usemann-{P}etty problem for arbitrary measures},
  Math. Ann. \textbf{331} (2005), no.~4, 867--887.

\bibitem{Zymonopoulou2008}
Marisa Zymonopoulou, \emph{The complex {B}usemann-{P}etty problem for arbitrary
  measures}, Arch. Math. (Basel) \textbf{91} (2008), no.~5, 436--449.

\bibitem{Zymonopoulou2009}
\bysame, \emph{The modified complex {B}usemann-{P}etty problem on sections of
  convex bodies}, Positivity \textbf{13} (2009), no.~4, 717--733.

\end{thebibliography}

\end{document}